\documentclass[11pt]{article}

\newcommand{\Z}{\mathbb{Z}}

\newcommand{\id}{|} 

\newcommand{\bradelta}{{\bf{ \Delta}} }
\newcommand{\bramu}{{\bf{m}} }

\usepackage{amsfonts,amssymb}
\usepackage{graphicx}

\newtheorem{theorem}{Theorem}[section]

\newtheorem{lemma}[theorem]{Lemma}
\newtheorem{proposition}[theorem]{Proposition}

\newtheorem{example}[theorem]{Example}

\newtheorem{remark}[theorem]{Remark}

\input{epsf.sty}

\begin{document}

\title{Algebraic Structures Derived from Foams}

\author{J. Scott Carter \\
University of South Alabama \\
\and
Masahico Saito\\
University of South Florida 
}

\maketitle

\begin{abstract}
Foams are surfaces with branch lines 
at which 
three sheets merge.  
They have been used in the 
categorification of $sl(3)$ quantum knot invariants and also  in physics.
The $2D$-TQFT  
of surfaces, on the other hand, is classified by means of 
commutative 
Frobenius algebras, 
where saddle points correspond to multiplication and comultiplication.
In this paper, we explore algebraic operations that branch lines derive 
under TQFT. In particular, we investigate 
Lie bracket and bialgebra structures. Relations to the original Frobenius algebra 
structures are discussed both algebraically and diagrammatically. 
\end{abstract}

\section{Introduction}

Frobenius algebras have been used extensively 
in the study of 
categorification of the Jones polynomial \cite{Kh06}, via $2$-dimensional 
Topological Quantum Field Theory ($2D$-TQFT, \cite{Kock}). 
For categorifications of other knot invariants, 
$2$-dimensional complexes called foams have been used instead \cite{Kh03,MV06}.
Although $2D$-TQFT has been characterized \cite{Kock} 
in terms of commutative Frobenius algebras, foams have not been algebraically characterized 
in terms of TQFT. Relations to Lie algebras, for example, have been 
suggested 
\cite{Kh03,MV06} through their boundaries 
which are 
called webs 
and 
that are trivalent graphs.
Foams have branch curves along which three sheets meet. Similar 
complexes 
appear as spines of $3$-manifolds, and 
 have been used for quantum invariants
\cite{CFS,ChFS,KSS,TV}.

Herein we study the 
types of algebraic operations 
that 
appear along the branch curves
of foams in relation to $2D$-TQFT. 
Recall that a $2D$-TQFT is a functor from the category of 
$2$-dimensional cobordisms to a category of $R$-modules 
(for some suitable ring $R$) that assigns an $R$-module to each connected 
component (circle) on the boundary of a surface, and an $R$-module homomorphism to a surface. 
In the case of a foam, 
we examine the associated
 algebraic operations that might be associated to branching circles in relation
  to the Frobenius algebra structure that occurs on the unbranched surfaces.
   Specifically, we identify 
and study Lie algebra and bialgebra structures in relation to branch curves,
and study their relations to 
the 
Frobenius algebra structure.

After reviewing necessary materials  in Section~\ref{prelimsec}, 
a Lie algebra structure 
for the branch curves  is  
studied in Section~\ref{Liesec}, and 
comultiplications  of bialgebras are 
examined 
 in Section~\ref{bialgsec}.
The foam skein theory based on 
the 
bialgebra case is also
defined 
in Section~\ref{bialgsec}.

\section{Preliminary}\label{prelimsec}

Algebraic structures we investigate  
include Frobenius algebras, Lie algebras and 
bialgebras. 
We restrict to the following situations.

A {\it Lie} algebra is a module $A$ over a unital commutative 
ring $R$ with a binary operation $[\ , \ ]: A \times A \rightarrow A$
that is bilinear, skew symmetric  ($[x, y]=-[y,x]$ for $x, y \in A$)
and satisfies the Jacobi identity
($[x,[y,z]]+[z,[x,y]]+[y,[z,x]]=0$ for $x,y,z\in A$). 

A {\it Frobenius algebra} is an algebra  over $R$ (that comes with associative linear 
multiplication $\mu: A \otimes A \rightarrow A$ and unit $\eta: R \rightarrow A$)
with  a non-degenerate 
form $\epsilon: A \rightarrow R$ that is associative ($\epsilon(x\otimes yz)=\epsilon(xy\otimes z)$
for $x,y,z \in A$). There is an induced co-associative comultiplication 
$\Delta: A \rightarrow A \otimes A$.
See \cite{CCEKS} for diagrams for Frobenius algebras
which we 
will 
use in this paper. 
A {\it bialgebra} is an algebra $A$ over $R$ with a comultiplication 
$\Delta: A \rightarrow A \otimes A$ that is an algebra homomorphism
($\Delta(xy)=\Delta(x)\Delta(y)$) and a counit $\epsilon: A \rightarrow R$
such that $(\epsilon \otimes {\rm id})\Delta= {\rm id} 
= ({\rm id}
\otimes \epsilon )\Delta$.
The following are typical examples.

\begin{example} \label{trancpolyex}
{\rm
Let $A=A_N$ be the Frobenius algebra of truncated polynomial 
$A_N=R [X] / (X^N)$ for a commutative unital ring $R$, with 
counit (Frobenius form) $\epsilon$  determined by 
$\epsilon (x^{N-1})=1$ and $\epsilon (x^i)=0$ for $i \ne N-1$.
The comultiplication $\Delta $ is determined by 
$\Delta(1) = \sum_{i=0}^{N-1} X^i \otimes X^{N-1-i}$. 
Diagrammatically, this is represented by a ``neck cutting'' relation~\cite{BN},
which we call 
{\it a $\Delta(1)$-relation} to distinguish the specific relation 
given  in \cite{BN} 
for $N=2$.    
See the right of Fig.~\ref{CNrel} for a diagrammatic representation of
the $\Delta(1)$-relation in this case. }\end{example}

\begin{figure}[htb]
\begin{center}
\includegraphics[width=1.5in]{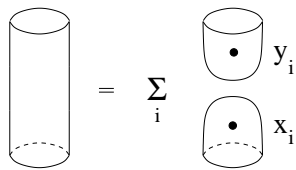}
\end{center}
\vspace*{-20pt}
\caption{The $\Delta(1)$-relation} 
\label{tube}
\end{figure}

In general, the $\Delta(1)$-relation is 
also described as follows (see \cite{Kh03,Kock}).  For  a commutative Frobenius algebra $A$ over a unital ring $R$ of finite rank
and with a non-degenerate Frobenius form $\epsilon$, 
there is a basis $\{ x_i\}$ and a dual basis $\{ y_i\}$, $i=1, \ldots, n$, 
such that $\epsilon (x_i y_i)=\delta_{i,j}$, the Kronecker   
delta, and $x=\sum_i y_i \epsilon (x_i x)$. 
This situation is depicted in Fig.~\ref{tube}, where the identity map $x \mapsto x$ 
in the LHS corresponds to the 
annular 
cobordism in the left of the figure, 
and the sum involving the Frobenius form $\epsilon$ is depicted in the right of the figure.

\begin{example} \label{MVex} {\rm
The Frobenius algebra structure on $A=\Z[a,b,c] [X] / (X^3-aX^2-bX-c)$ 
is presented  
in \cite{MV06} as follows.
The multiplication and the unit are defined by those for polynomials, 
the Frobenius form (counit)  $\epsilon$ is defined by 
$\epsilon(1)=\epsilon(X)=0$, $\epsilon(X^2)=-1$. 
The comultiplication is accordingly computed as 
\begin{eqnarray*}
\Delta (1) &=& -( 1 \otimes X^2 + X \otimes X + X^2 \otimes 1 )
+ a ( 1 \otimes X + X \otimes 1) + b (1 \otimes 1),  \\
\Delta (X) &=& - (X \otimes X^2 + X^2  \otimes X) + a (X \otimes X) - c (1 \otimes 1),  \\
\Delta (X^2) &=& - (X^2 \otimes X^2) - b( X\otimes X) -c (1 \otimes X+ X\otimes 1) .
\end{eqnarray*}
}
\end{example}

\begin{figure}[htb]
\begin{center}
\includegraphics[width=2in]{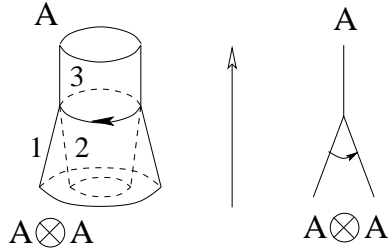}
\end{center}
\vspace*{-20pt}
\caption{Operation on a branch circle}
\label{branch}
\end{figure}

We fix a $2D$-TQFT such that a connected circle corresponds to $A$. 
For TQFTs we refer to \cite{Kock}. 
We follow definitions of foams in \cite{Kh03,MV06}, 
except that facets of foams are decorated by basis 
elements of $A$, in a general way as in \cite{Kaiser}. 
A foam without boundary is called closed.

\begin{figure}[htb]
\begin{center}
\includegraphics[width=.8in]{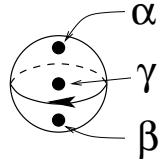}
\end{center}
\vspace*{-20pt}
\caption{The theta foam}
\label{theta}
\end{figure}

We briefly summarize their definitions.
${\bf Foam}_A$  is the category of formal linear combination over $R$ of
cobordisms 
of compact $2$-dimensional complexes 
in $3$-space
with the following data.
(1) Boundaries are planar graphs with trivalent rigid vertices.
(2) For an interior point of a cobordism,  the  neighborhood of each point  is homeomorphic 
to either Euclidean 
$2$-space (a facet)  or a {\it branch  curve} where 
three facets  of half planes meet.
(3) Each facet 
is oriented, and the induced orientation 
on the branch curve is consistent among three facets that share the curve.
(4) A cyclic order of facets are specified using the orientation of  
$3$-space
as depicted in Fig.~\ref{branch}. 
(5) Each facet has a basis element of $A$ assigned.
(6) 
An annular cobordism
as depicted in the left of  Fig.~\ref{tube} 
is equivalent to the linear combination as depicted 
on 
the right. 
(7) Values $\theta(\alpha, \beta, \gamma) \in A$ 
of the theta foam, 
as depicted in Fig.~\ref{theta} are specified.

In \cite{Kh03,MV06}, 
it was shown that the values in $A$ of closed foams are well-defined 
for values of the theta foams, as long as the 
{\it cyclic symmetry condition} 
$\theta(\alpha, \beta, \gamma)=\theta( \beta, \gamma, \alpha)=\theta(\gamma,\alpha, \beta)$
is satisfied.

By $2D$-TQFT for a chosen $A$,
the 
 two circles  
 on 
the left Fig.~\ref{branch}
are mapped to 
the factors of 
$A \otimes A$. For the cyclic order along the oriented branch circle as depicted, 
make a correspondence between the facet labeled $1$, $2$, $3$, respectively, 
to the first, second, and the target factor of $A \otimes A \rightarrow A$. 
Thus the cobordism near a branch circle as depicted in the figure induces 
a linear map $A \otimes A \rightarrow A$ under the chosen TQFT and 
the values of theta foams.  
Denote this map by $\bramu: A \otimes A \rightarrow A$. 
The
goal of this paper 
is 
to investigate this map. 

In terms of maps among tensor products of $A$s, we use 
planar graphs regularly used in knot theory, as well as Frobenius 
algebras as in \cite{CCEKS}. 
In particular, the Frobenius form (the counit) is depicted by a maximum, 
unit by a minimum, (co)multiplications by trivalent vertices. 
In this convention, diagrams are read from bottom to top, corresponding to 
the domain and range of maps. 
The map $\bramu$ corresponding to theta foams has 
a 
specified 
cyclic order,  
as indicated 
on the right of 
Fig.~\ref{branch}.  
The map $\bramu$ is defined with this specific 
order, and the map with the opposite order, depicted by a diagram 
with the opposite arrow, represent the map $\bramu \circ \tau$, 
where 
$\tau: A \otimes A \rightarrow A \otimes A$ is the map induced from the transposition
$\tau(x \otimes y) = y \otimes x$.

\section{Lie  algebras}\label{Liesec}

In this section we 
show that there are infinitely many TQFTs under which Lie algebra structures 
are induced from 
the branch circle 
operation.
Since our goal is to exhibit a Lie bracket, 
in this section we use the notation $[\ , \ ]  : A \times A \rightarrow A$,
instead 
of $\bramu: A \otimes A \rightarrow A$.

\begin{proposition}
For any commutative unital  ring $R$ and a positive odd integer $N > 1$,  
there exist a Frobenius algebra $A$ over $R$ and
values of the theta foams in ${\bf Foam}_A$
such that the branch circle operation $\bramu$ induces a non-trivial Lie algebra structure on  $A$.
 \end{proposition}
{\it Proof.\/}
Let $A=R[X]/(X^N)$ for an odd integer $N>1$. 
For simplicity we denote $\theta (X^a, X^b, X^c)$ by $\theta (a,b,c)$ in this proof.

Let $N>3$. 
Define $\theta( a,b,c)=1$ if $a=0$, $b+c=N$ and $1<b<c$, 
as well as  all cyclic permutations of such $(a,b,c)$. 
Define $\theta( a,b,c) =-1$ if $a=0$, $b+c=N$ and $1<c<b$, 
as well as  all cyclic permutations of such $(a,b,c)$. 
Finally define $\theta( a,b,c)= 0$ for all the other cases.
For $N=3$, replace the conditions $1<b<c$ and $1<c<b$, 
respectively, by $b<c$ and $c<b$. 
 We show that these theta foam values induce Lie brackets 
 as desired.

\begin{figure}[htb]
\begin{center}
\includegraphics[width=3in]{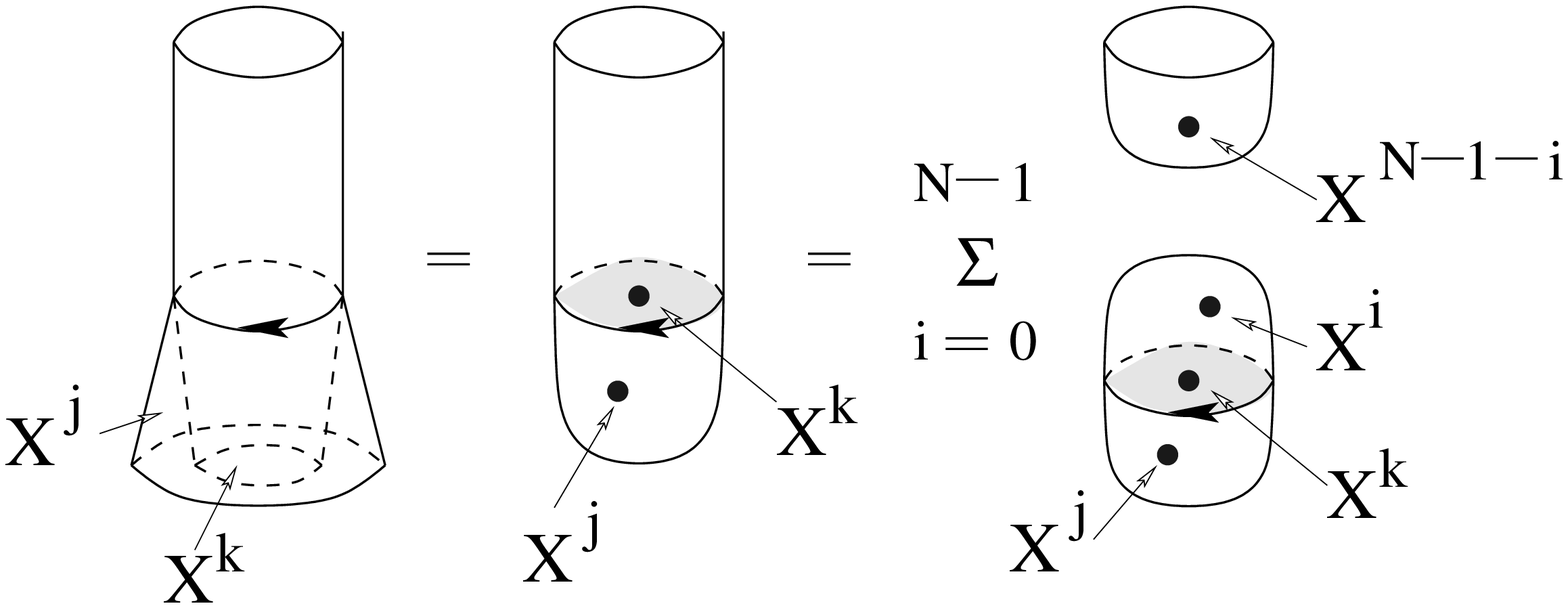}
\end{center}
\vspace*{-20pt}
\caption{Evaluating bracket}
\label{CNrel}
\end{figure}

The operation $[X^j, X^k]$ is evaluated, using the $\Delta(1)$-relation,  by 
$$[X^j, X^k]= \sum_{i=0}^{N-1} \theta (X^i, X^j, X^k) X^{N-1-i}. $$
This calculation is depicted 
in Fig.~\ref{CNrel}.
Since $\theta (i,j,k) =0$ unless $i+j+k=N$, we have 
$[X^j, X^k]=  \theta (i,j,k) X^{N-1-i}$ where $i=N-(j+k)$ and $N-1-i=j+k-1$, so that 
$[X^j, X^k]=  \theta ({N-(j+k)}, j, k) X^{j+k-1}$. Note that if $j+k>N$, then 
the RHS is understood to be zero from the definition of $\theta$. 
{}From the definition of $\theta$ by cyclic ordering, the skew symmetry 
of $[\ , \ ]$ is clear.  We show the  Jacobi identity 
$$ [ X^j, [ X^k, X^\ell ] ] +  [ X^\ell, [ X^j, X^k]] +  [ X^k, [ X^\ell, X^j]] =0$$
case by case. First we compute
\begin{eqnarray*}
 [ X^j, [ X^k, X^{\ell } ]] 
 &=& \theta ({N -(k+\ell )}, k, \ell) \ \theta( N+1-(j+k+\ell), j, k +\ell -1) , \\ {}
   [ X^\ell, [ X^j, X^k]] &=&  \theta ({N -(j+k )}, 
   j, k)\  \theta( N+1-(j+k+\ell), \ell, j+k-1),   \\ {}
    [ X^k, [ X^\ell, X^j]] &=&  \theta ({N -(\ell+j )}, 
\ell , j)\  \theta( N+1-(j+k+\ell), k,  \ell +j -1), 
    \end{eqnarray*}
hence it is sufficient to prove that the sum of the right-hand sides is zero.

\noindent
{\bf Case 1}: $j+k+\ell > N+1$. 

In this case, the second 
factors  
of the RHS are  zero, so that all terms are zero.

\noindent
{\bf Case 2}: $j+k+\ell \leq  N+1$ and $k + \ell > N$. 

This case implies that $j=0$ and $k+\ell=N+1$. Since $N+1$ is even, 
$k$ and $\ell$ have the same parity. The first 
factor $\theta ({N -(k+\ell)}, k, \ell) $ is $0$ 
since 
$N -(k+\ell)=-1$.  
(When the arguments of $\theta$ are out of range, then $\theta=0$. )
Suppose $k=\ell=(N+1)/2$.  
Then the second and the third  terms are 
\begin{eqnarray*}
\lefteqn{
\theta ({N -(j+k)}, 
j, k)\  \theta( N+1-(j+k+\ell), \ell, j+k-1) } \\
&=&  \theta ((N-1)/2, 0, (N+1)/2 )\  \theta (0,  (N+1)/2,  (N-1)/2 ) = (-1) (-1) = 1, \\
\lefteqn{
 \theta ({N -(\ell+j)},  
 \ell , j)\  \theta( N+1-(j+k+\ell), k,  \ell +j -1) } \\
&=&   \theta ((N-1)/2,  (N+1)/2, 0 )\  \theta (0,  (N+1)/2,  (N-1)/2 ) = (1) (-1)=-1,
  \end{eqnarray*}
  as desired. Hence assume $k < \ell$ without loss of generality.
For the second and third terms, we have 
\begin{eqnarray*}
\lefteqn{
 \theta ({N -(j+k)} 
 j, k) \ \theta( N+1-(j+k+\ell), \ell, j+k-1) } \\
&=&  \theta ( \ell -1 , 0, k) \ \theta (0, \ell, k-1) = (1) (-1) = -1, \\
\lefteqn{
 \theta ({N -(\ell+j)}, 
  \ell , j)\  \theta( N+1-(j+k+\ell), k,  \ell +j -1) }\\
&=&  \theta ( k -1 , \ell, 0)\  \theta (0, k, \ell-1) = (1) (1) = 1,
  \end{eqnarray*}
  as desired. 

\noindent
{\bf Case 3}: $j+k+\ell \leq  N+1$ and $k + \ell \leq  N$. 

First we check the case where two of $j, k, \ell$ are the same. 
Suppose $j=k$. Then the second term is zero, as 
$\theta(N-(j+k), j, k)=0$. 
Furthermore, for the first and third terms, we have 
$\theta ({N -(k+\ell)}, k, \ell) = -  \theta ({N -(k+\ell)}, \ell , j)$ and 
$$ \theta( N+1-(j+k+\ell), j, k +\ell -1)=
 \theta( N+1-(j+k+\ell), k,  \ell +j -1)$$  
  as desired.
 The other cases $(k=\ell, j=\ell)$ are checked similarly.
 Hence we may assume $j < k < \ell$. 
 
 Since $\theta$ vanishes unless one of the entries is $0$, 
the first  
factors 
of the RHS are zero if $N> k+ \ell$ and  $0 < j<k<\ell$. 
 hence we may assume that $k+\ell = N$ or $j=0$.

We continue to examine specific subcases. Suppose that $j=0$ and $k + \ell < N$.  The RHS becomes:
\begin{eqnarray*}\lefteqn{
\theta(N-(k+\ell), k, \ell)\theta(N+1-(k+\ell),0,k+\ell-1)} \\
&+& \theta(N-k, 0, k)\theta(N+1-(k+\ell),\ell,k-1) \\
& +& \theta(N-\ell, \ell,0)\theta(N+1-(k+\ell),k,\ell-1)  .
 \end{eqnarray*}
 
 If $j=0$ and $k+\ell = N$, we have 
$$\theta(0, k, \ell)\theta(1,0,N-1) +
\theta(\ell, 0, k)\theta(1,\ell,k-1) +
\theta(k, \ell,0)\theta(1,k,\ell-1)  . $$

If $1<k$, then the sum is $0$ since $\theta(1,0,N-1)=0$ and the arguments of 
the second and the third factors are all non-zero. If $k=1$, then $\theta(0,1,N-1)=\theta(\ell,0,1)=\theta(1, \ell,0)$, so the sum is $0$.

Now suppose that $k+\ell<N$. The first term is $0$ and the second factors in the sum have arguments that are all non-zero unless $k=1$. If $k=1$, we have
$$\theta(N-1,0,1)\theta(N-\ell,\ell, 0)+ \theta (N-\ell, \ell, 0)\theta(N-\ell,1,\ell-1)=0 . $$

Finally, suppose that $j\ne0$, so that $k+\ell=N$. The RHS becomes:
$$\theta(0,N-\ell, \ell)\theta(1-j,j, N-1) + \theta(0,j,N-\ell)\theta(1-j,\ell, N-1) + \theta(N-(\ell +j),\ell, j)\theta(1-j,N-\ell,\ell+j-1).$$
If $1<j$, then the first argument of all the second factors is negative, so the sum is $0$. If $j=1$, then each term has a factor that is $0$. 
$\Box$

 \bigskip

Since the original motivation came from the foams in \cite{Kh03,MV06}, 
we examine the 
Frobenius algebra   
in \cite{MV06} closely. 
In this case, the multiplication that is induced by branch circles also satisfies the Jacobi identity.

\begin{proposition}
\label{fromVaz}
Let $A=\Z[a,b,c] [X] / (X^3-aX^2-bX-c)$ 
with Frobenius structure defined as  in Example~\ref{MVex} from \cite{MV06}.
The branch curve operation 
 $[ \ , \ ] $ is skew-symmetric and 
satisfies the 
Jacobi identity:
$$ 
[ U , [ V, W ] ]
+  
[ V , [ W, U] ]
 +  
 [ W, [ U, V ] ]
 $$
 for any $U,V,W \in A$. 
\end{proposition}
{\it Proof.\/} 
This is confirmed by calculations. {}From the axioms of $A$ and the theta foam values 
 that are given in \cite{MV06}:  
$$\theta(1,X,X^2)=\theta(X^2,1,X)=\theta(X,X^2,1)=1=-\theta(1,X^2,X)
=-\theta(X,1,X^2)=
-  
\theta(X^2,X,1)$$
 while $\theta=0$ for any other arguments,  
we compute using the $\Delta(1)$ relation for Example~\ref{MVex}:

\begin{eqnarray*}
[1, X ] &=& -1 ,  \\ {}
[1, X^2 ] &=& X- a , \\ {}
[X, X^2] &=& -X^2 + aX + b . 
\end{eqnarray*}
Then one computes 
\begin{eqnarray*}
[ 1, [X, X^2]] &=& -X , \\ {}  
[X, [X^2, 1]] &=& a , \\ {}
[X^2, [1, X]] &=& X-a ,
\end{eqnarray*}
as desired.  
In general, we consider cyclic permutations of $X^j,$ $X^k$, and $X^\ell$ in the expression $[X^j ,[X^k, X^\ell]]$. Since the bracket is skew-symmetric, then we need only consider the cases in which $j$, $k$, and $\ell$ are distinct. The remaining case follows by skew-symmetry.
$\Box$

\begin{figure}[htb]
\begin{center}
\includegraphics[width=2in]{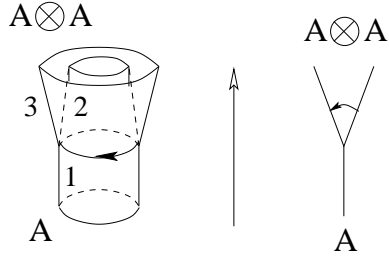}
\end{center}
\vspace*{-20pt}
\caption{Upside-down operation}
\label{branchup}
\end{figure}

\begin{figure}[htb]
\begin{center}
\includegraphics[width=1.7in]{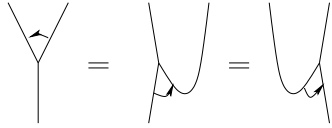}
\end{center}
\vspace*{-20pt}
\caption{$\bradelta$ can be defined from left or right}
\label{leftright}
\end{figure}

We 
define 
the operation 
$\bradelta:  A \rightarrow A \otimes A $ 
that is associated to the left of Fig.~\ref{branchup}, 
 a diagram that is  up-side down of Fig.~\ref{branch},
 in which one circle branches into two from bottom to top.
 A cyclic order is specified in the figure. If we specify the ordered tensor factors
 assigned to each sheet by $A_i$, $i=1,2,3$, 
 then the operation is defined as 
 $\bradelta:  A_1 \rightarrow A_3 \otimes A_2 $. 
 A planar diagram representing this operation is depicted in the right of the figure. 
Imitating  Sweedler
notation 
$\Delta(u)=\sum u_{(1)} \otimes u_{(2)} $
for comultiplication, 
we denote $\bradelta (u)=\sum u_{ ((1))} \otimes u_{((2))}.$
The next lemma relates this operation to the unit map, 
 and diagrammatic formulations are given in Fig.~\ref{leftright}.

\begin{lemma} 
Let $A=\Z[a,b,c] [X] / (X^3-aX^2-bX-c)$, with $\Delta(1)$-condition defined as in Example~\ref{MVex}. 
The map $\bradelta$ is computed as follows.
\begin{eqnarray*}
\bradelta (1) &=& 1 \otimes X - X \otimes 1 , \\
\bradelta (X) &=& a( 1  \otimes X - X \otimes 1 ) - (1  \otimes X^2  - X^2 \otimes 1 ),  \\ 
\bradelta (X^2) &=& (a^2 + b) ( 1  \otimes X - X \otimes 1 ) - a  (1  \otimes X^2  - X^2 \otimes 1 ) 
+ (X  \otimes X^2  - X^2 \otimes X ) . 
\end{eqnarray*}
\end{lemma}

Direct calculations show

\begin{lemma}
$\bradelta (V ) = \sum [ V, 1_{(1)} ] \otimes 1_{(2)}= \sum 1_{(1)} \otimes [ 1_{(2)}, V ]$.
\end{lemma}

\begin{figure}[htb]
\begin{center}
\includegraphics[width=4in]{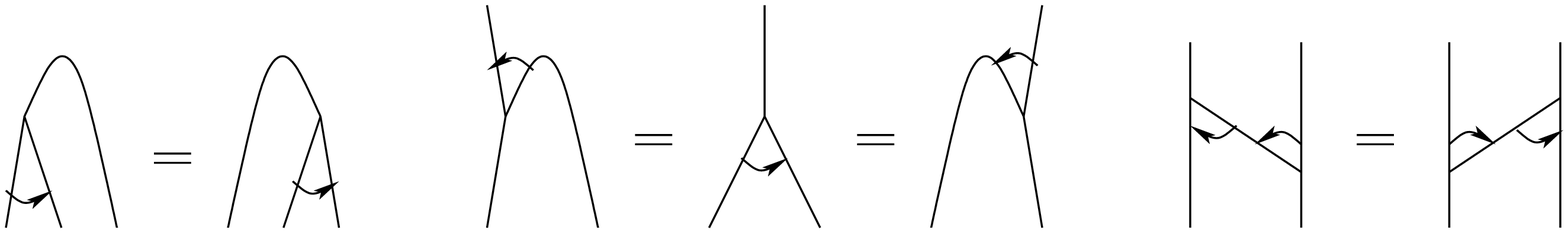}
\end{center}
\vspace*{-20pt}
\caption{Other symmetric relations}
\label{otherrels}
\end{figure}

The diagram for this relation is depicted in Fig.~\ref{leftright}. Other relations
that follow are depicted in Fig.~\ref{otherrels}.

\begin{figure}[htb]
\begin{center}
\includegraphics[width=5in]{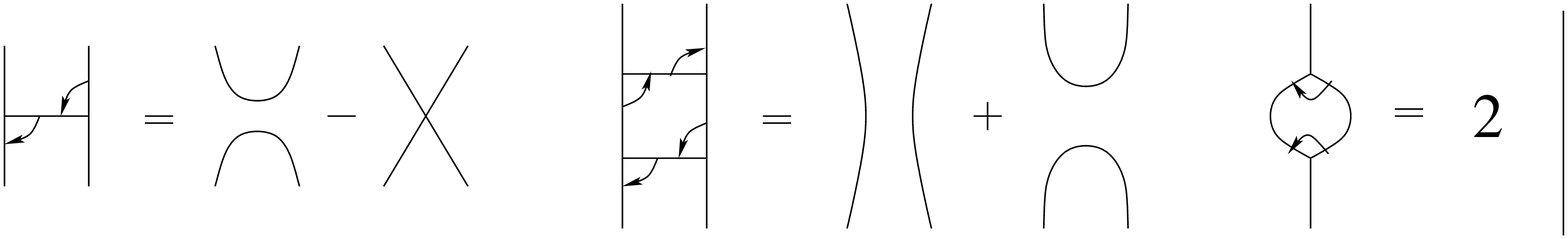}
\end{center}
\vspace*{-20pt}
\caption{Web skein relations}
\label{webskeins}
\end{figure}

The following relations hold for maps in Frobenius algebras and maps associated to branch circles.
Here we used the notation $\bramu$ instead of $[\ , \ ]$ to formulate in tensor products.
The equalities are diagrammatically represented in Fig.~\ref{webskeins}.

\begin{figure}[htb]
\begin{center}
\includegraphics[width=3.5in]{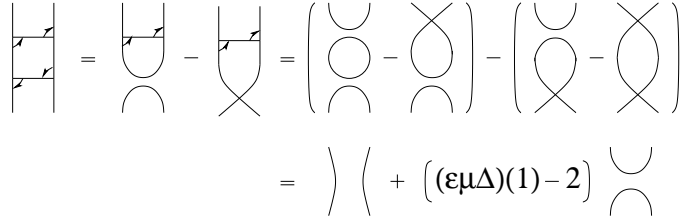}
\end{center}
\vspace*{-20pt}
\caption{Proof of the skein relation}
\label{webskeinproof}
\end{figure}

\begin{proposition}\label{skeinprop}
For $A=\Z[a,b,c]/(X^3-aX^2-bX-c)$ with $\theta$ values as above, 
the map $\bradelta : A \rightarrow A \otimes A$ satisfies the following identities:
\begin{eqnarray*}
(\bramu \otimes \id )(\id \otimes \bradelta) &=& 
 \Delta(1) ( \epsilon \mu )-  \tau  , \\ 
 ( \ (\bramu \otimes \id )(\id \otimes \bradelta) \ ) ^2 &=& \id + \Delta(1) ( \epsilon \mu ), \\ 
\bramu \bradelta &=& 2 \ \id . 
\end{eqnarray*}
\end{proposition}
{\it Proof.\/} 
The first  
and the third equalities are verified 
by calculations on basis elements.
For all $X^i$ and $X^j$, it is computed as
$[X^i, X^j_{((1))}] \otimes X^j_{((2))} = X^j \otimes X^i + \epsilon (X^{i+j} ) \Delta (1)$. 
The second 
relation is 
diagrammatically computed as in Fig.~\ref{webskeinproof}.
Note that the {\it handle element}  
 $\epsilon\mu\Delta(1)$ is $3$. 
$\Box$

\begin{remark} {\rm
The skein relations stated in Proposition~\ref{skeinprop},
as planar diagrams (instead of surface skein relation),  
coincide with those described in 
\cite{Kh03} as a description of Kuperberg's invariant \cite{Kup91}, 
with the choice of $q=1$.

Thus, the operation at branch curve of the foam used to categorify the 
quantum $sl(3)$ invariant satisfies the skein relations at the classical limit 
of the invariant. 
}
\end{remark}

\begin{figure}[htb]
\begin{center}
\includegraphics[width=4in]{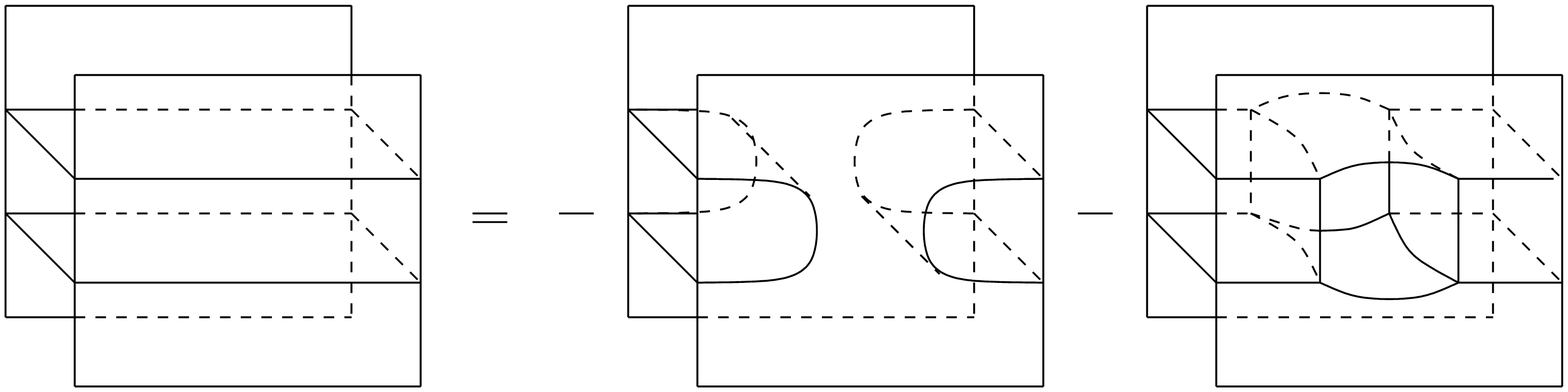}
\end{center}
\vspace*{-20pt}
\caption{A surface skein relation in \cite{Kh03,MV06} }
\label{sqr}
\end{figure}

\begin{remark} {\rm
The second relation in Proposition~\ref{skeinprop} is related to 
the 
local surface 
skein relation in \cite{Kh03,MV06} as follows. 
Their local relation is depicted in Fig.~\ref{sqr}. 
Notice the negative signs, as well as resemblance to our relation.
After performing their relations locally, move the holes 
of each term along the $S^1$ factor to the other side.
Then one obtains a tube connecting two sheets.
Then perform the {\it bamboo cutting relation}, 
that is computed by applying $\Delta(1)$-relation three times. 
In this case, one computes that it is the negative 
of the original bamboo segment. 
These negative signs cancel, and we obtain our equation. 
Thus, our relation follows from theirs, or algebraically  as we have shown.

We also point out that the second and the third relation in
Proposition~\ref{skeinprop} have interpretations 
 in ${\bf Foam}_A$.
One simply takes the product of these diagrams with $S^1$ to obtain 
foams, and the equalities hold in ${\bf Foam}_A$.
The first equality, however, is not realized in ${\bf Foam}_A$, 
as the intersection of surfaces are not allowed in ${\bf Foam}_A$.
} \end{remark}

\begin{figure}[htb]
\begin{center}
\includegraphics[width=4in]{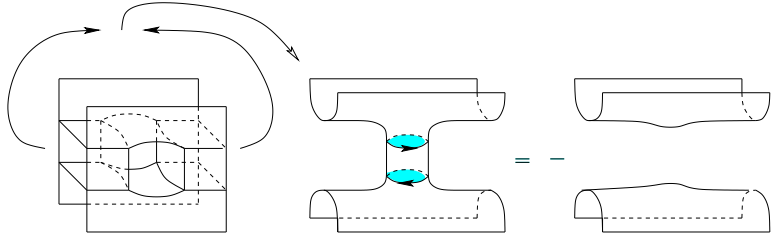}
\end{center}
\vspace*{-20pt}
\caption{Cutting a bamboo segment}
\label{bambocut}
\end{figure}

\section{Bialgebras}\label{bialgsec}

In this section, we investigate functors whose image of branch curves
induce bialgebra structure 
for 
group algebras. 
Let $G$ be a group.
Let $A=R[G]$ be the group ring with a commutative unital ring $R$. 
It is well known that $A$ has a commutative 
Hopf algebra structure defined as follows (see, for example, \cite{Kock}).
Define $\bradelta: A \rightarrow A \otimes A$ by linearly extending 
$\bradelta(x)=x \otimes x$. 
(This is different from the comultiplication as a Frobenius algebra 
$\Delta(x)= \sum_{x=yz}  y \otimes z$.)
The unit map is defined as the same as the Frobenius unit map 
$\eta (1)=1_G$, where $1_G$ is the identity element of $G$.
(The counit map as a Frobenius algebra is defined by $\epsilon (1_G)=1$ and  $\epsilon(x)=0$
for $x \neq 1_G$.)
The following shows that there is a strong requirement for group algebras to give 
bialgebra structures through branch curves.

\begin{figure}[htb]
\begin{center}
\includegraphics[width=2.2in]{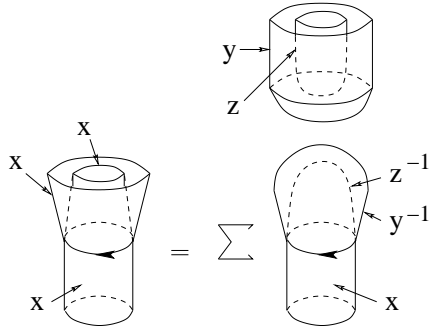}
\end{center}
\vspace*{-20pt}
\caption{Comultiplication by theta foams}
\label{bialg}
\end{figure}

\begin{proposition}\label{bialgprop}
Let $G$ be an abelian group.
For any 
unital 
ring $R$, 
 the branch circle operation $\bramu$ induces  a bialgebra structure on $A$
if and only if every non-identity element of $G$ has order $2$. 
 \end{proposition}
{\it Proof.\/}
The $\Delta(1)$-relation is written as $\Delta(1)=\sum_{y \in G} y \otimes y^{-1}$, 
and the reducing $\bradelta$ into the theta foam is depicted in Fig.~\ref{bialg}. 
For $\bradelta (x)=x \otimes x$ to hold in the figure, we have $y=z=x$, 
and the value of the theta foam 
being 
$\theta(x, y^{-1}, z^{-1})=1$ for 
$y=z=x$ and $0$ otherwise. 

For $\theta$ to satisfy the cyclic symmetry, 
this condition is satisfied if and only if $x^{-1}=x$
($x$ having order $2$) for any $x\in G$, 
and in this case, the theta foam values are determined by 
$\theta (x,x,x)=1$ for any $x \in G$ and $0$ otherwise.
$\Box$

\begin{figure}[htb]
\begin{center}
\includegraphics[width=2.2in]{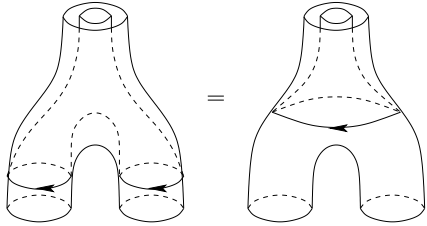}
\end{center}
\vspace*{-20pt}
\caption{The compatibility condition of a bialgebra}
\label{compatisfce}
\end{figure}

\begin{remark} {\rm
The condition of a bialgebra that the comultiplication is 
an algebra homomorphism
(also called a compatibility condition)
 $\bradelta (a b)=\bradelta (a) \bradelta (b)$ 
for $a, b \in A$, is represented 
by  
surfaces in Fig.~\ref{compatisfce}. 
}
\end{remark}

\begin{figure}[htb]
\begin{center}
\includegraphics[width=3.3in]{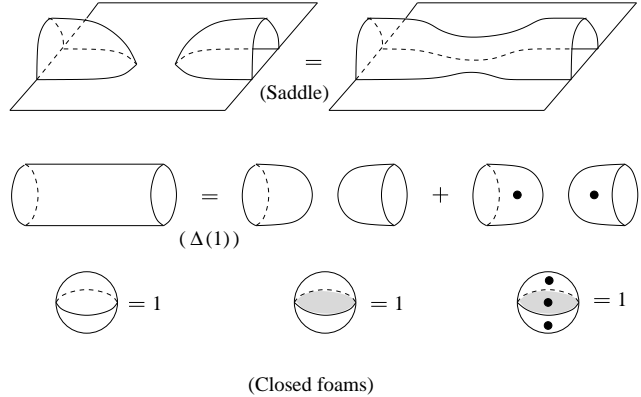}
\end{center}
\vspace*{-20pt}
\caption{Surface skein relations}
\label{skein}
\end{figure}

\begin{remark}{\rm
In \cite{AF}, skein modules for $3$-manifolds based on embedded surfaces 
modulo the surface skein relations described in \cite{BN} 
were defined and studied. 
Surface skein modules were generalized in \cite{Kaiser} using 
general commutative Frobenius algebras.
Such notions are directly generalized to foams, with various skein relations at hand. 
Skein modules for $sl(3)$ foams are analogously defined using 
the local skein relations given in \cite{Kh03,MV06}, for example.

Here we propose  local skein relations based on the foams in Proposition~\ref{bialgprop}
with  the bialgebra on branch  
curves for  the group ring $\Z[x]/(x^2-1)$. 
Considering
that 
 the move characteristic to bialgebras is the compatibility condition as 
depicted in Fig.~\ref{compatisfce}, 
we take a local change that happens at the saddle point of this move, 
as depicted in the top of Fig.~\ref{skein} (labeled as saddle) 
as a local surface skein relation. Other relations in Fig.~\ref{skein}
are those coming from Frobenius algebra structure and the theta foam values as before.

Thus 
 the skein module $ {\bf F}(M) $  
in this case  can be defined  
to be the isotopy classes of foams in a given $3$-manifold $M$
modulo the local surface skein relations in Fig.~\ref{skein}. 
Although computations of this skein module in general 
 is out of the scope of this paper,
it seems 
interesting to look 
into relations between foams and 
the  
 topology of $3$-manifolds.

}\end{remark}

\subsection*{
Acknowledgments}

JSC was supported in part by NSF Grant DMS \#0603926. MS was supported in
part by NSF Grants DMS \#0603876 and \#0900671.


\begin{thebibliography}{99}
\setlength{\itemsep}{-5pt}

\bibitem{AF}
M.\,Asaeda; C.\,Frohman, 
{\it A note on the Bar-Natan skein module}, 
arXiv:math/0602262.

\bibitem{BN} D.\,Bar-Natan,  
 {\it Khovanov's homology for tangles and cobordisms}, 
 Geom. Topol. {\bf 9} (2005) 1443--1499.
 
 \bibitem{CCEKS}
  J.S.\,Carter,  A.S.\,Crans, M.\,Elhamdadi, E.\,Karadayi, M.\,Saito, 
{\it  Cohomology of Frobenius Algebras and the Yang-Baxter Equation,}
 Communications in Contemporary Mathematics
 {\bf 10}, Suppl. 1 (2008),  791--814.

\bibitem{CFS}
J.S.\,Carter, D.\,Flath, M.\,Saito,
 {\it  Classical and Quantum 6j Symbols,} 
Mathematical notes, vol. 43, 
 Princeton University Press,  1995.



\bibitem{ChFS} S.\,Chung, M.\,Fukuma, A.\,Shapere,
{\it Structures of topological lattice field theories in three dimensions,} 
Internat. J. Modern Phys. A  {\bf 9}  (1994),   1305--1360.

\bibitem{Kaiser}
U.\,Kaiser,
{\it Frobenius algebras and skein modules of surfaces in $3$-manifolds},
{arXiv:0802.4068}.


\bibitem{KSS}
L.H.\,Kauffman, M.\,Saito, M.C.\,Sullivan,
{\it Quantum invariants of templates,}
J. Knot Theory Ramifications, {\bf 12} (2003), 653-681.


\bibitem{Kh99} 
M.\,Khovanov, 
{\it A categorification of the Jones polynomial},
{Duke Math. J.},
{\bf 101(3)} (1999), 
359--426.


\bibitem{Kh03}
M.\,Khovanov, 
{\it $sl(3)$ link homology},
Algebraic \& Geometric Topology, 
{\bf 4} (2004),  1045--1081.

\bibitem{Kh06}
M.\,Khovanov, 
{\it Link homology and Frobenius extensions},
{Fundamenta Mathematicae}, 
{\bf 190}  (2006), 179--190.

\bibitem{Kock} J.\,Kock,
{\it Frobenius algebras and 2D topological quantum field theories,}
London Mathematical Society Student Texts {\bf 59}, Cambridge University Press, 2003.

\bibitem{Kup91}
G.\,Kuperberg,
{\it Spiders for rank $2$ Lie algebras,} Comm. Math. Phys., {\bf 180(1)}  (1996), 109--151.


\bibitem{MV06}
M.\,Mackaay; P.\,Vaz, 
{\it The universal $sl_3$-link homology},
arXiv:math/0603307.


\bibitem{TV}  V.\,Turaev, O.\,Viro,  {\it State Sum Invariants of 
3-Manifolds and Quantum 6J-Symbols}, Topology {\bf 31}, No 4 (1992),
865--902.


\end{thebibliography}
\end{document}